\documentclass{amsart}

\usepackage{amsfonts}       
\usepackage{amsthm}
\usepackage{amssymb}
\usepackage{xcolor}
\usepackage{soul}
\usepackage[pdftex]{graphicx}
\usepackage{caption}
\usepackage{float}

\newtheorem{theorem}{Theorem}[section]

\theoremstyle{definition}

\theoremstyle{remark}

\numberwithin{equation}{section}

\newcommand{\N}{{\mathbb{N}}}

\newcommand{\R}{{\mathbb{R}}}

\newcommand{\cico}{C^{\infty}_c(\Omega)} 
\newcommand{\cic}{C^{\infty}_c}

\newcommand{\darr}[4]{{\left\{\begin{array}{ll}
   {#1}&{#2}\\[0.1cm]
   {#3}&{#4}
 \end{array}\right.}}

\newcommand{\be}{\begin{equation}}
\newcommand{\ee}{\end{equation}}
\newcommand{\bea}{\begin{eqnarray}}
\newcommand{\eea}{\end{eqnarray}}
\newcommand{\bean}{\begin{eqnarray*}}
\newcommand{\eean}{\end{eqnarray*}}
\newcommand{\la}{\label}

\newcommand{\lphc}{  \Big(\frac{p-1}{p}\Big)^p}

\begin{document}

\title{The Hardy constant: a review}

\author{G. Barbatis}
\address{Department of Mathematics, National and Kapodistrian University of Athens}
\email{gbarbatis@math.uoa.gr}






\begin{abstract}
We present a review of
results that have been obtained in the past
twenty-five  years concerning the $L^p$-Hardy inequality with 
distance to the boundary. We concentrate on results 
where the best Hardy constant is either computed exactly or estimated from below.
\end{abstract}

\maketitle

\section{Introduction}

In this artricle we present a review of some 
of the results that have been obtained in the past twenty-
five years concerning the Hardy inequality
\be
\la{hi1}
\int_{\Omega}|\nabla u|^pdx \geq c \int_{\Omega}
\frac{|u|^p}{d^p}\, dx \; , \qquad u\in C^{\infty}_{c}(\Omega).
\ee
Here $\Omega$ is an open and connected subset of $\R^n$ 
with non-empty boundary and $d(x)={\rm dist}(x,\partial\Omega)$, 
$x\in\Omega$, denotes the distance to the boundary of $\Omega$.
Hardy inequalities involving the function $d(x)$ are sometimes
called geometric Hardy inequalities in order to distinguish them
from Hardy inequalities involving the distance to an interior
point.

Since the publication of the review article \cite{D99} 
the literature related to inequality \eqref{hi1} has grown significantly.
Several aspects of this inequality
as well as other related
inequalities have been extensively studied: weighted 
inequalities,
Rellich inequalities, improved inequalities,
inequalities on non-Euclidean settings, fractional inequalities 
and more. The  publication of three books \cite{BEL15,GM13,Ru19} 
is indicative of the recent interest on this area.

In the present article we shall be primarily concerned with
the best constant for inequality \eqref{hi1}. Hence we
shall present results where  the best constant is 
precisely computed as well as results where lower 
estimates are obtained. Some mention of improved Hardy
inequalities will also be made. At the end of the article we present
some open problems.

\section{The Hardy constant}

The $L^p$-Hardy inequality involving the distance to the boundary reads
\be
\la{hi}
\int_{\Omega}|\nabla u|^pdx \geq c \int_{\Omega}
\frac{|u|^p}{d^p}\, dx \; , \qquad u\in C^{\infty}_{c}(\Omega).
\ee
Here $p>1$, $\Omega\subset\R^n$ is a domain and $d(x)=
{\rm dist}(x,\partial\Omega)$, $x\in\Omega$.

We say that the $L^p$-Hardy inequality is valid for the domain
$\Omega$ if there exists $c>0$ such that \eqref{hi}
holds true. We denote by $H_p(\Omega)$ the best 
constant for \eqref{hi}, the $L^p$-Hardy 
constant of the domain $\Omega$. In case $p=2$ we
shall simply write $H(\Omega)$.

There are various sufficient conditions as well as 
necessary conditions for the validity of the Hardy 
inequality. These are typically related to some 
regularity of the domain. 
As already mentioned, we shall be primarily concerned
with the precise value of the Hardy 
constant as well as with explicit lower estimates.

\subsection{Domains with critical Hardy constant}

For $p>1$ we set
\[
\alpha_p = \lphc .
\]
This constant plays a special role for the $L^p$-Hardy 
inequality. Besides being the Hardy constant in dimension 
one, it is also the case that
$H_p(\Omega)\leq  \alpha_p$
when some part of $\partial\Omega$ is
$C^2$; indeed a lower regularity is enough, see
\cite[Theorem 5]{MMP98}.

The importance of the value $\alpha_p$
is also indicated by the following dichotomy which
has been obtained in \cite{MMP98,MaSh00}.
\begin{theorem}
Let $\Omega$ be a bounded domain with $C^2$ boundary.
There holds $H_p(\Omega)<\alpha_p$ if and only if
the $L^p$ Hardy quotient admits a minimizer in
$W^{1,p}_0(\Omega)$.
\end{theorem}
This has been generalized in \cite{LP19} to bounded domains with $C^{1,\gamma}$ boundary. We note that for the `only if'
part lower boundary regularity is enough, see \cite{Te98}.

There are several conditions under which the $L^p$
Hardy constant is equal to $\alpha_p$. The one that is
probably best known is the convexity of the domain
$\Omega$. This has been long known in case $p=2$ while
the case of general $p>1$ was obtained in \cite{MaSo97}.

A more general condition was established in \cite{BFT04} where
it was proved that if  the domain $\Omega$ is such that
\be
\Delta d \leq 0 , \quad \mbox{ in }\Omega
\la{wmc}
\ee
(in the distributional sense) then
\[
 \int_{\Omega}| \nabla u |^p dx \geq \lphc
\int_{\Omega}\frac{|u|^p}{d^p}dx
\; , \quad u\in C^{\infty}_{c}(\Omega).
\]
Indeed the following series improvement
was established in \cite{BFT03a}. Define
$X_1(t) =(1-\log t)^{-1}$, $X_{k}(t)=X_1(X_{k-1}(t)$,
$k\geq 2$, $t\in (0,1)$. (These are iterated
logarithmic functions that vanish
at $t=0$ at a rate that becomes slower as $k$ increases.)

\begin{theorem}
\la{thm:series}
Assume that for the domain $\Omega\subset\R^n$
condition \eqref{wmc} is satisfied and also assume
that $\sup_{\Omega} d<+\infty$. Then for any $p>1$ there exists a constant $D\geq \sup_{\Omega}d$ such that
\[
 \int_{\Omega}| \nabla u |^p dx \geq \lphc
\int_{\Omega}\frac{|u|^p}{d^p}dx +
\frac{p-1}{2p} \Big(\frac{p-1}{p}\Big)^{p-2}
\sum_{k=1}^{\infty}\int_{\Omega}
\frac{|u|^p}{d^p}X_1^2X_2^2\ldots X_k^2dx,
\]
for all $u\in C^{\infty}_{c}(\Omega)$; here
$X_j=X_j(d(x)/D)$. Moreover the inequality
is sharp as each new term of the series is added.
\end{theorem}
Theorem \ref{thm:series} is one amongst several results where an 
improvement of a Hardy inequality with sharp constant is 
obtained. To our knowledge the first such result
is contained in \cite{Ma85}. Improved Hardy inequalities can be 
of two types: the added term may contain a weighted $L^p$ norm or 
some  weighted $L^q$ norm with $p<q\leq p^*$
where $p^* =np/(n-p)$ is the Sobolev exponent. Improvements of 
the first type, also called homogeneous improvements,
have been obtained, among others, 
in \cite{BEL12,BFT04,BFT03a,BFT03b,BM97,FTT09,HO2L02,T05}.
Concerning Sobolev improvements see
\cite{BFT04,BFL08,BM97,FMaT07,FTT09,FL12,Gk13};
we shall not go into any further details regarding improved Hardy inequalities.

Going back to condition \eqref{wmc}, we note that it is also mentioned in 
\cite{D99}. Domains for which \eqref{wmc} is valid are called 
weakly mean convex domains. 
Any convex domain is weakly mean convex.
It has been proved independently in \cite{LLL12,Ps13}
that if $\Omega$ is bounded with
$C^2$ boundary then weak mean convexity is equivalent to 
mean convexity, that is to the mean curvature of
$\partial\Omega$ being non-negative.
In two dimensions and for bounded domains with $C^2$
boundary mean  convexity is equivalent to convexity \cite{AK85}.
This is not true in higher dimensions.

In case $p=2$ there are other domains that are known to have Hardy constant equal to $1/4$.
One such domain is the annulus
$\{ x\in\R^n  : r<|x| <R \}$ in dimension $n\geq 3$ \cite{MMP98}.
This covers also the case $R=\infty$, i.e.
the complement of a ball.
The latter case can also be obtained as a special case of a more general
result in \cite{Gk13} where it is shown that if the domain $\Omega\subset\R^n$
satisfies
\[
 -\Delta d + (n-1) \frac{\nabla d\cdot x}{|x|^2} \geq 0 \; , \qquad 
 \mbox{ in }\Omega
\]
then the Hardy inequality with constant $1/4$ is valid for $\Omega$.

Using different methods the following
theorem has been obtained in \cite{Av13,Av14}:
\begin{theorem}
Let $n\geq 2$. There exists a number $\Lambda_n > 0$ (expressed by means
of a certain hypergeometric equation) such that:
if $\Omega\subset\R^n$ is a bounded domain and if for 
each point $y\in\partial\Omega$ there exists a ball 
$B_r$ of radius $r=(\sup_{\Omega}d)/\Lambda_n$ 
with $y\in \partial B_r$
and $B_r \subset \Omega^c$, then $H(\Omega) =1/4$.
\end{theorem}
An analogous result for $p\in (1,2)$ is contained in \cite{Av15}.

The following theorem about cones has been proved in \cite{D95}; see also \cite{DPP17}.
\begin{theorem}
\la{thm:cone}
Let $U$ be an open connected subset of the unit sphere and let $\Omega$ be the corresponding infinite cone, given in spherical coordinates by
$\Omega=\{(r,\omega) : r>0 \, , \; \omega\in U\}$.
Let $\delta:U\to\R$ be such that $d=r\delta(\omega)$ on $\Omega$. Then the Hardy constant of $\Omega$ coincides with the best constant $k_U$ for the inequality
\[
\int_U |\nabla_{\omega}g|^2dS +
\Big( \frac{n-2}{2} \Big)^2 \int_U 
g^2dS  \geq k_U \int_U \frac{g^2}{\delta^2}dS ,
\qquad g\in \cic(U).
\]
\end{theorem}

\section{The Hardy constant in two dimensions}

There is more that can be said about the $L^2$-Hardy constant 
when we consider domains in $\R^2$. This is due to the availability of tools from complex 
analysis but also to the fact that 
explicit computations are simpler than in higher 
dimensions.

The following well known result of Ancona \cite{Anc86}
is proved by a suitable application of Koebe's 1/4 
theorem.
\begin{theorem}
If $\Omega\subset\R^2$ is a simply connected domain then 
\[
\int_{\Omega}|\nabla u|^2 dx \geq \frac{1}{16} 
\int_{\Omega} \frac{u^2}{d^2}\, dx 
\; ,\qquad u\in \cico .
\]
\end{theorem}
In \cite{LS08} a modified version of Koebe's 1/4 theorem
was used to prove the next theorem which involves a
quantified measure of non-convexity.
\begin{theorem}
Let $\Omega\subset\R^2$ be simply connected and satisfy an external cone condition: each $y\in\partial\Omega$ is the vertex of an infinite cone $C$ of angle $\theta$ with $\Omega\subset C$. Then
\[
\int_{\Omega}|\nabla u|^2dx \geq \frac{\pi^2}{4\theta^2} \int_{\Omega} \frac{u^2}{d^2}\, dx \; , \quad \mbox{ for all }u\in C^{\infty}_c(\Omega).
\]
\end{theorem}
We note that if $\Omega$ is convex then one recaptures the constant 1/4.

\subsection{Sectors in $\R^2$}

Let $\beta \in [\pi, 2\pi]$ and
let $\Lambda_{\beta}$ denote the infinite sector of angle $\beta$,
\[
 \Lambda_{\beta} =\{ (r,\theta) : r>0 \, , \;\; 0<\theta<\beta\}.
\]
Then $d^{-2}=r^{-2}V_{\beta}(\theta)$ where
\[
V_{\beta}(\theta) =
\left\{ \begin{array}{lll} \frac{1}{\sin^2\theta}, & 0< \theta  < \frac{\pi}{2} ,  \\
1 , & \frac{\pi}{2}<\theta <\beta-\frac{\pi}{2} , \\
\frac{1}{\sin^2(\beta-\theta)}, & \beta-\frac{\pi}{2} <\theta <\beta. 
\end{array}
\right.
\]
It follows from Theorem \ref{thm:cone} that
the Hardy constant $H(\Lambda_{\beta})$ coincides
with the best constant $c_{\beta}$ for the Hardy-type inequality
\[
\int_0^{\beta} g'(\theta)^2 d\theta \geq
c_{\beta}
\int_0^{\beta} g(\theta)^2 V_{\beta}(\theta) d\theta \; , 
\quad  g\in \cic(0,\beta).
\]
Equivalently, $c_{\beta}$ is the largest constant $c$ for which
the boundary value problem
\be
\la{bvp}
 \darr{ \psi''(\theta) +c V_{\beta}(\theta)\psi(\theta) =0  , \;\;\;}{ 0\leq\theta\leq\beta,}
{ \psi(0)=\psi(\beta)=0,}{}
\ee
has a positive solution in $(0,\beta)$.

The qualitative behavior of $c_{\beta}$ as a function of $\beta$ was studied in 
\cite{D95} where it was proved that there exists a critical angle $\beta_{cr}\in (4,2\pi)$ such that 
$c_{\beta}=1/4$ if $\beta\in [\pi,\beta_{cr}]$ while $c_{\beta}$ is strictly 
decreasing in the interval $[\beta_{cr},2\pi]$.
Numerical computations give
$\beta_{cr} \simeq 1.546\pi$ and $c_{2\pi}\simeq 0.205$.
In \cite{BT14} the boundary value problem \eqref{bvp}
was further analyzed to obtain the following explicit
description of $c_{\beta}$:
\begin{theorem} The critical angle $\beta_{cr}$ is the unique 
solution in  the interval $(\pi,2\pi)$ of the equation
\[ 
\tan\big( \frac{\beta_{cr}-\pi}{4}\big) = 4\bigg(\frac{\Gamma(\frac{3}{4})}{\Gamma(\frac{1}{4})}\bigg)^2.
\]
Moreover for  $\beta\in [\beta_{cr},2\pi]$ the 
constant $c_{\beta}$ is the unique solution to the 
equation
\[
\sqrt{c_{\beta}}\tan\big( \sqrt{c_{\beta}}(\frac{\beta-\pi}{2})\big) = 2 \bigg(
\frac{\Gamma(\frac{3+\sqrt{1-4c_{\beta}}}{4})}
{\Gamma(\frac{1+\sqrt{1-4c_{\beta}}}{4})}
\bigg)^2 .
\]
\end{theorem}
The next theorem was obtained in \cite{BT14}
and provides an explicit description of the
Hardy constant of an arbitrary non-convex quadrilateral.
\begin{theorem}
\la{thm:quad}
The Hardy constant of a non-convex quadrilateral is equal
to $c_{\beta}$ where $\beta\in (\pi,2\pi)$ is the
size of the non-convex angle.
\end{theorem}
The proof of this theorem makes a combined use of the distance function and the solution $\psi(\theta)$ of \eqref{bvp}. Following this method the Hardy constants 
of other domains were computed in \cite{BT15}.

\newpage

\noindent
{\bf Some open problems.} 

\vspace{0.2cm}

We have made a brief exposition of some of the results that have been obtained in 
the past twenty-five years concerning geometric Hardy inequalities.
Naturally, several interesting problems remain open.
We close this short review article by presenting three 
such problems.

\medskip

\noindent
{\em Problem 1.} Let
$H^*$ be the largest constant
such that
\[
\int_{\Omega}|\nabla u|^2 dx \geq H^* \int_{\Omega} \frac{u^2}{d^2}\, dx \
\]
for all simply connected domains $\Omega\subset\R^2$
and for all $u\in\cico$. Determining the exact value of $H^*$ 
is an open problem first posed, to our knowledge, in \cite{L05};
see also \cite{Ban99}. 
What is currently known is that
$H^* \in [\frac{1}{16} , c_{2\pi}]$ and the exact computation of $H^*$
seems very difficult.

A more realistic immediate target would be to determine 
the Hardy constant of more planar domains and/or to 
narrow the above interval for $H^*$. The fact that for the proof
of Theorem \ref{thm:quad} five different types of 
quadrilaterals were distinguished is indicative 
of how challenging even this simpler problem is.

\medskip

\noindent
{\em Problem 2.} Consider the inequality
\be
\la{whc}
\int_{\Omega}|\nabla u|^2 dx \geq \alpha
\int_{\Omega} \frac{u^2}{d^2} -\beta \int_{\Omega} 
u^2 dx
\ee
where $\alpha,\beta$ are real numbers and $u\in \cico$.
The number
\[
H_w(\Omega) =\sup \big\{ \alpha\in\R  :
\mbox{there exists $\beta\in\R$ 
such that \eqref{whc} holds for all $u\in\cico$}
\big\} 
\]
is called the weak Hardy constant of the domain
$\Omega\subset\R^n$.

The weak Hardy constant has been studied in
detail in \cite{D95}. It was shown in particular
that for any bounded domain $\Omega$ there holds
$
H_w(\Omega) =\min\{ h(y) \; : \; y\in\partial\Omega\}
$
where $h(y)$ is the local Hardy constant at the point
$y\in\partial\Omega$ defined as
\bean
h(y) &=&\lim_{r\to 0} \; \sup \big\{ \alpha\in\R \, :
\mbox{there exists $\beta\in\R$ 
such that }\\
&&\hspace{3cm}\mbox{\eqref{whc} holds for all $u\in\cico \cap
B_r(y)$}\big\} .
\eean
It was conjectured in \cite{D95} that $h(y)\leq 1/4$
for any bounded domain $\Omega$ and any $y\in\partial\Omega$. This remains open to this day. Indeed we are not aware of any result stating that
$H(\Omega)\leq 1/4$ for any bounded domain $\Omega\subset\R^n$, a weaker version of the above conjecture.

\medskip

\noindent
{\em Problem 3.} From \cite[Theorem 1]{BT06} and
\cite[Exercise 4.2.10]{Gr14} easily follows the following
\begin{theorem}
Let $p>1$ and $n,m\in\N$, $n,m\geq 2$. There exists a constant
$c(n,m,p)$ such that for any weakly mean convex domain $\Omega\subset\R^n$ and
any $u\in \cico$ there holds
\footnote{Here $|\Delta^{m/2} u|$ stands for
$|\nabla\Delta^{(m-1)/2} u|$ when $m$ is odd.
}
\be
\int_{\Omega}|\Delta^{m/2} u|^p dx \geq c(n,m,p)\int_{\Omega}
\frac{|u|^p}{d^{mp}}dx 
\la{lp:rel}
\ee
\end{theorem}
The best value of the constant $c(m,n,p)$ is not known.
Even for the seemingly simple inequality
\[
\int_{\R^n_+}|\Delta u|^p dx \geq c\int_{\R^n_+}
\frac{|u|^p}{x_n^{2p}}dx  \; , \qquad u\in\cic(\R^n_+)
\]
the best constant is not known unless $p= 2$.
This is in sharp contrast to the case of Rellich inequalities with distance to an interior point where for $m\in\N$ and $p>1$ satisfying $mp<n$ the best
constant $A(m,p,n)$ for the inequality
\[
\int_{\Omega}| \Delta^{m/2} u|^p dx \geq A(m,p,n)
\int_{\Omega} \frac{|u|^p}{|x|^{mp}} dx \; , 
\quad \cico ,
\]
is known \cite{DH98} and an infinite series 
improvement with iterated logarithms similar to that
of Theorem \ref{thm:series} has been obtained 
in the framework of a general Cartan-Hadamard 
manifold \cite{Ba07}.
We note that in case $p=2$ and for convex $\Omega$ the best constant
in \eqref{lp:rel} and a sharp series improvement have been 
obtained in \cite{Ba06,BT06,O99}. Whether the latter 
estimates are valid for mean convex domains is also an open problem.

\

\noindent
{\bf Acknowledgement.} I thank A. Tertikas for helpful comments.

\


\end{document}